\documentclass[a4paper]{classe3-1}

\usepackage{amsmath,amssymb}
\usepackage{pstricks,pst-node}
\usepackage[francais,english]{babel}
\usepackage[T1]{fontenc}
\usepackage[latin1]{inputenc}

\newtheorem{theorem}{Theorem}[section]
\newtheorem{lemma}[theorem]{Lemma}
\newtheorem{e-proposition}[theorem]{Proposition}

\newtheorem{e-definition}[theorem]{Definition\rm}
\newtheorem{remark}{\it Remark\/}
\newtheorem{example}{\it Example\/}

\newtheorem{theoremebis}{Th\'eor\`eme}

\renewcommand{\theequation}{\arabic{equation}}
\def\og{\leavevmode\raise.3ex\hbox{$\scriptscriptstyle\langle\!\langle$~}}
\def\fg{\leavevmode\raise.3ex\hbox{~$\!\scriptscriptstyle\,\rangle\!\rangle$}}

\begin{document}
\sffamily
\begin{frontmatter}


\selectlanguage{english}
\title{The fundamental group of a triangular algebra without double bypasses.}
\vspace{-2.6cm}
\selectlanguage{francais}
\title{Le groupe fondamental d'une alg\`ebre triangulaire sans double raccourci.}
\selectlanguage{english}
\author[authorlabel1]{Patrick Le Meur}
\ead{lemeur@math.univ-montp2.fr}
\address[authorlabel1]{Institut de Math\'ematiques et de Mod\'elisation de Montpellier, UMR CNRS 5149, 
Universit\'e Montpellier 2, case courier 051, place Eugène Bataillon, 34095 Montpellier cedex 5, France}
\begin{abstract}
%
%
%
%
%
%
Let $A$ be a basic connected finite dimensional algebra over a field of characteristic zero.
A fundamental group depending on the presentation of $A$ has been defined by several authors (see \cite{martinezvilla}).
Assuming the quiver of $A$ has no oriented cycles, no double arrows and no double by-passes, 
we show there exists a suitable presentation of $A$ with quiver and admissible relations, with fundamental group
denoted by $\pi_1(A)$, such that the fundamental
group of any other presentation of $A$ with quiver and admissible relations is a quotient of $\pi_1(A)$. 
{\it To cite this article: P. Le Meur, C. R.
Acad. Sci. Paris, Ser. I ??? (200?).}
\vskip 0.5\baselineskip
\selectlanguage{francais}
%
%
%
%
%
%
\noindent{\bf R\'esum\'e}
\vskip 0.5\baselineskip
\noindent
Soit $A$ une alg\`ebre basique connexe et de dimension finie sur un corps de caract\'eristique nulle.
Plusieurs auteurs (voir \cite{martinezvilla}) ont d\'efini pour $A$ un groupe fondamental d\'ependant du choix d'une 
pr\'esentation de $A$.
En supposant que le carquois de $A$ n'a pas de cycle orient\'e, pas de fl\`eches
parall\`eles et pas de double raccourci, nous d\'emontrons qu'il existe une pr\'esentation privil\'egi\'ee
de $A$ par carquois et relations admissibles, de groupe fondamental not\'e $\pi_1(A)$,
 telle que le groupe fondamental de toute autre pr\'esentation de $A$ par carquois et relations admissibles est un quotient
de $\pi_1(A)$. {\it Pour citer cet article~: P. Le Meur, C. R.
Acad. Sci. Paris, Ser. I ??? (200?).}
\end{abstract}
\end{frontmatter}
%
%
%
%
%
%
\selectlanguage{francais}
\section*{Version fran\c{c}aise abr\'eg\'ee}
Soit $A$ une alg\`ebre basique, connexe et de dimension finie sur un corps $k$,
et soit $Q$ le carquois (connexe) de $A$.
Si $I$ est un idéal admissible de $kQ$ tel que $A\simeq kQ/I$,
le groupe fondamental $\pi_1(Q,I)$ de $(Q,I)$ a \'et\'e d\'efini
dans \cite{martinezvilla} à partir de la relation d'homotopie $\sim_I$ de $(Q,I)$.
Il existe des exemples d'alg\`ebres $A$ avec différentes pr\'esentations admissibles $A\simeq kQ/I$ et $A\simeq kQ/J$
telles que $\pi_1(Q,I)\not\simeq\pi_1(Q,J)$. Le but de ce texte est d'étudier cette situation 
et la possible existence d'un groupe fondamental qui serait canoniquement attaché à cette algèbre.
Pour cela nous considérons les raccourcis de $Q$ c'est à dire les couples $(\alpha,u)$ o\`u $\alpha$ 
est une fl\`eche de $Q$ et $u$ est
un chemin (orient\'e) de $Q$ parall\`ele \`a $\alpha$ et distinct de $\alpha$. Un double raccourci de $Q$ est un quadruplet
$(\alpha,u,\beta,v)$ o\`u $(\alpha,u)$ et $(\beta,v)$ sont deux raccourcis tels que $\beta$ est une fl\`eche
parcourue par le chemin $u$. Dans cette note, nous montrons le th\'eor\`eme suivant:
\begin{theoremebis}
Supposons que le corps $k$ est de caractéristique nulle. Supposons que $Q$ n'a pas de cycle orient\'e et que $Q$ n'admet pas de
double raccourci. Alors il existe une pr\'esentation admissible $kQ/I_0\simeq A$ de $A$
telle que pour toute autre pr\'esentation admissible $kQ/I\simeq A$ il existe un morphisme
surjectif de groupes $\pi_1(Q,I_0)\twoheadrightarrow\pi_1(Q,I)$.
\end{theoremebis}
Notons que l'hypoth\`ese sur les doubles raccourcis implique que le carquois $Q$ n'admet pas de
fl\`eches parall\`eles.
Si $kQ/I\simeq A$ et $kQ/J\simeq A$ sont deux pr\'esentations admissibles de $A$ telles que $\sim_I$ est plus fine que $\sim_J$
(i.e. $\gamma\sim_I\gamma'\Rightarrow\gamma\sim_J\gamma'$), alors il existe un morphisme surjectif de groupes
$\pi_1(Q,I)\twoheadrightarrow\pi_1(Q,J)$. Partant de cette remarque, la preuve du Th\'eor\`eme consiste \`a construire un carquois
$\Gamma$ dont les sommets sont les relations d'homotopie des pr\'esentations admissibles de $A$ et tel que si $\sim_I\to\sim_J$
est une fl\`eche de $\Gamma$, alors $\sim_I$ est strictement plus fine que $\sim_J$. Avant de d\'ecrire les fl\`eches de $\Gamma$, nous introduisons
des automorphismes particuliers de $kQ$: les transvections et les dilatations. Une dilatation est un automorphisme 
$\varphi\colon kQ\xrightarrow{\sim}kQ$ tel que $\varphi(\alpha)\in k^*\alpha$ pour toute fl\`eche $\alpha$. Une transvection est un
automorphisme de la forme $\varphi_{\alpha,u,\tau}\colon kQ\xrightarrow{\sim}kQ$ o\`u $(\alpha,u)$ est un raccourci, $\tau\in k$ et
$\varphi_{\alpha,u,\tau}$ est d\'efini par $\varphi_{\alpha,u,\tau}(\alpha)=\alpha+\tau u$ et $\varphi_{\alpha,u,\tau}(\beta)=\beta$ pour toute fl\`eche $\beta\neq\alpha$.
L'int\'er\^et des dilatations et des transvections est le suivant: soit $kQ/I\simeq A$ une pr\'esentation admissible de $A$,
soit $\varphi\colon kQ\to kQ$ un automorphisme et soit $J=\varphi(I)$. Si $\varphi$ est une dilatation alors $\sim_I$ et $\sim_J$ co\"incident,
 si $\varphi$ est une transvection alors l'une des deux relations d'homotopie $\sim_I$ ou $\sim_J$ est plus fine (au sens large) que l'autre.
Cette propri\'et\'e permet de d\'efinir les fl\`eches de  $\Gamma$: il existe une fl\`eche $\sim\to\sim'$
dans $\Gamma$ si et seulement si il existe deux pr\'esentations admissibles $kQ/I\simeq A$ et $kQ/J\simeq A$
ainsi qu'une transvection $\varphi\colon kQ\to kQ$ telles que $\sim=\sim_I$, $\sim'=\sim_J$, $J=\varphi(I)$ et $\sim_I$ est strictement plus 
fine que $\sim_J$. De cette fa\c{c}on, $\Gamma$ est un carquois connexe, la longueur des chemins orient\'es de $\Gamma$
est born\'ee et tout sommet de $\Gamma$ est le but d'un chemin orient\'e dont la source est une source de $\Gamma$ (i.e. n'est le but
d'aucune fl\`eche de $\Gamma$). Nous montrons alors
 que le carquois $\Gamma$ n'a qu'une seule source et que si $kQ/I_0\simeq A$ est une pr\'esentation admissible de $A$
telle que $\sim_{I_0}$ est l'unique source de $\Gamma$, alors le couple $(Q,I_0)$ satisfait la conclusion du Th\'eor\`eme.
Il est \`a noter que si $car(k)\neq 0$, il existe des exemples d'alg\`ebre $A$ dont le carquois $Q$ n'a ni cycle orient\'e ni
 double raccourci et telle que le carquois $\Gamma$ admet plusieurs sources (le Th\'eor\`eme reste cependant vrai pour ces exemples).
La preuve de ce Th\'eor\`eme sera d\'etaill\'ee dans un article en pr\'eparation qui fera partie de la th\`ese de l'auteur
\`a Montpellier sous la direction de Claude Cibils, le cadre de travail sera alors \'elargi au contexe des 
rev\^etements galoisiens de l'alg\`ebre $A$.
%
%
%
%
%
%
\selectlanguage{english}
%
%
%
%
%
%
\section{Introduction}
Let $A$ be a finite dimensional connected algebra over a field $k$.
We are interested in the representation theory of $A$, thus we may assume $A$
is basic and we denote by $Q$ its ordinary quiver.
Any presentation $kQ/I\simeq A$ of $A$ with quiver and admissible relations
(i.e. $I$ is an admissible ideal of $kQ$ and $kQ/I\simeq A$ is an isomorphism of $k$-algebras)
gives rise to a group ($\pi_1(Q,I)$) called the fundamental group of $(Q,I)$ 
(see \cite{bongartz} and \cite{martinezvilla}).
The fundamental group is particularly useful in covering techniques (see
 \cite{bongartz},\cite{green1} and \cite{green2}). However
different presentations of $A$ with quiver and admissible relations may provide non isomorphic fundamental groups.
In this text we intend to clarify this \textit{uncanonical} situation. 
For this purpose we will use the notion of double bypass. Recall a bypass (see \cite{assem}) is a pair
$(\alpha,u)$ where $\alpha$ is an arrow of $Q$ and $u$ is an oriented path of $Q$
parallel to $\alpha$ and distinct from $\alpha$. A double bypass is quadruple $(\alpha,u,\beta,v)$
where $(\alpha,u)$ and $(\beta,v)$ are bypasses such that the arrow $\beta$ appears in the path $u$.
Assume that
$car(k)=0$ and suppose the algebra $A$ satisfies:
\begin{enumerate}
\item[.] $A$ is triangular (i.e. $Q$ has no oriented cycles),
\item[.] the quiver $Q$ has no double bypasses.
\end{enumerate}
\begin{theorem}
\label{THM}Assuming the above conditions,
there exists a presentation $kQ/I_0\simeq A$ with quiver and relations such that for any other
presentation $kQ/J\simeq A$ there is a surjective group morphism $\pi_1(Q,I_0)\to\pi_1(Q,J)$.	
\end{theorem}
\vskip 10pt
Notice that our assumption on double bypasses implies that $Q$ has no double arrows
(i.e. $A$ is square free). Indeed, if $\alpha\neq\beta$ are parallel arrows, then $(\alpha,\beta,\beta,\alpha)$
is a double bypass. Notice also that if we assume moreover
that $A$ is shurian, we recover the result \cite[Thm. 3.5]{marcos}.\\
This note is part of the author's thesis made at Université Montpellier 2 under the supervision of
 Claude Cibils. A detailed version of this note will be written in a subsequent paper and in the author's 
thesis report.
%
%
%
%
%
%
\section{Basic definitions and notations}
Let $(Q,I)$ be a connected \textbf{quiver with admissible relations}, i.e. $Q$ is a (finite) connected quiver $Q$ and $I$ is an admissible ideal of $kQ$. Recall that admissible means
$(kQ^+)^n\subseteq I\subseteq (kQ^+)^2$ for some $n$ (where $(kQ^+)^n$ stands for the ideal generated by 
paths of length $n$).
A \textbf{walk} is an unoriented path of $Q$. The stationnary walk at a vertex $x$ will be denoted by $e_x$.
Let $r=t_1u_1+\ldots+t_nu_n\in I$ where $t_i\in k^*$ and the $u_i$'s are distinct paths, then
$r$ is called a \textbf{minimal relation} if $n\geqslant 2$ and if for any 
non empty proper subset $E$ of $\{1,\ldots,n\}$
the term $\sum_{i\in E}t_i.u_i$ does not lie in $I$.
The \textbf{homotopy relation} $\sim_I$ of $(Q,I)$ is the smallest equivalence relation on the set of walks (of $Q$)
 which is compatible with the concatenation of walks and such that:
\begin{enumerate}
\item[.] for any arrow $\alpha$ with source $x$ and target $y$, 
we have $\alpha\alpha^{-1}\sim_I e_y$ and $\alpha^{-1}\alpha\sim_I e_x$.
\item[.] $u_1\sim_I u_2$ as soon as $t_1u_1+\ldots+t_nu_n$ is a minimal relation of $(Q,I)$.
\end{enumerate}
Let then $x_0$ be a vertex of $Q$. 
The \textbf{fundamental group} (see \cite{bongartz} and \cite{martinezvilla}) $\pi_1(Q,I,x_0)$ of $(Q,I)$ at $x_0$ is the set
of $\sim_I$-homotopy classes of walks starting and ending at $x_0$. The composition in $\pi_1(Q,I,x_0)$ is induced
by the concatenation of walks and the unit is the $\sim_I$-class of $e_{x_0}$. As the choice of 
$x_0$ is unrelevant (since $Q$ is connected) we will write $\pi_1(Q,I)$ for short.
\begin{example}
Assume $Q$ is equal to
\psset{unit=10pt}
\begin{pspicture}(0,0)(6,3)
\rput(0,0){\rnode{A}{}}
\rput(2,2){\rnode{B}{}}
\rput(4,0){\rnode{C}{}}
\rput(6,0){\rnode{D}{}}
\ncline{->}{A}{B}
\aput{0}{$b$}
\ncline{->}{A}{C}
\bput{0}{$a$}
\ncline{->}{B}{C}
\aput{0}{$c$}
\ncline{->}{C}{D}
\aput{0}{$d$}
\end{pspicture}
and set $I=<\ da\ >$ and $J=<\ da-dcb\ >$.
\vskip 10pt
 Then $kQ/I\simeq kQ/J$ whereas $\pi_1(Q,I)\simeq\mathbb{Z}$ and $\pi_1(Q,J)= 0$.
\end{example}
%
%
%
%
%
%
\section{Preliminary results}
As we wish to compare the fundamental group of different presentations of a given algebra,
it is natural to try to compare the corresponding homotopy relations. In this section 
we give two helpful lemmas for such a comparison. Let us first introduce some terminology.
Let $Q$ be any quiver with set of vertices $Q_0$ and set of arrows $Q_1$.
Let $\varphi\colon kQ\to kQ$ be an automorphism which is the identity map on $Q_0$. We will say that
$\varphi$ is a \textbf{dilatation} if $\varphi(\alpha)\in k^*\alpha$ for any $\alpha\in Q_1$.
We will say that $\varphi$ is a \textbf{transvection} if there exists 
a bypass $(\alpha,u)$ and $\tau \in k$ such that $\varphi=\varphi_{\alpha,u,\tau}$ where 
$\varphi_{\alpha,u,\tau}(\alpha)=\alpha+\tau u$ and $\varphi_{\alpha,u,\tau}(\beta)=\beta$ 
for $\beta\in Q_1\backslash \{\alpha\}$.
In  analogy with the classical decomposition 
of elements of $GL_n(k)$ as products of transvections and dilatations, 
we have the following Lemma.\\
\begin{lemma}
\label{lem1}
Assume $Q$ has no oriented cycles and let $(Q,I)$ and $(Q,J)$ be two quivers with admissible relations
satisfying $kQ/I\simeq kQ/J$.
Then there exists $\psi\colon kQ/I\to kQ/J$ an isomorphism which is the identity map on $Q_0$.
Moreover for any such $\psi$, there exists $\varphi\colon kQ\to kQ$ an automorphism such that 
$\varphi(I)=J$ and such that $\varphi$ induces $\psi$.
Finally $\varphi$ is a composition of a dilatation and finitely many transvections.
\end{lemma}
\vskip 10pt
Note that Lemma \ref{lem1} can (partially) be rewritten in terms of semi-direct product of groups
as follows: the group of automorphisms $kQ\to kQ$ which are the identity map on $Q_0$
is the semi-direct product of the subgroup of dilatations with the (normal) subgroup generated
by transvections. Notice also that other studies
of the automorphism group of an algebra where made in relation with the Picard group
of the algebra (see \cite{guil-asensio}, \cite{pollack} and \cite{strametz}).\\
We now turn to a fundamental Lemma which is the first important step in the 
proof of Theorem \ref{THM}.\\
\begin{lemma}
\label{lem2}
Assume $Q$ has no oriented cycles and let $(Q,I)$ and $(Q,J)$ be quivers with admissible relations.\\
Let $\varphi\colon kQ\to kQ$ be an automorphism with $\varphi(I)=J$.
If $\varphi$ is a dilatation then $\sim_I$ and $\sim_J$ coincide.
Assume now $\varphi=\varphi_{\alpha,u,\tau}$ is a transvection.
\begin{enumerate}
\item[a)] If $\alpha\sim_I u$ and $\alpha\sim_J u$ then $\sim_I$ and $\sim_J$ coincide.
\item[b)] If $\alpha\not\sim_I u$ and $\alpha\sim_J u$ then $\sim_J$ is generated by $\sim_I$ and $\alpha\sim_J u$.
\item[c)] If $\alpha\not\sim_I u$ and $\alpha\not\sim_J u$ then $\sim_I$ and $\sim_J$ coincide and $I=J$.
\end{enumerate}
\end{lemma}
\vskip 10pt
\begin{remark}
The word \textbf{generated} stands for: genarated as an equivalence relation
 which is compatible with the concatenation of walks an such that $\alpha\alpha^{-1}\sim e_y$
and $\alpha^{-1}\alpha\sim e_x$ for any arrow $x\xrightarrow{\alpha}y\in Q_1$.
\end{remark}
\vskip 5pt
\begin{remark}
The following implication (symmetrical to b)):
\begin{center}
{
\textit{If $\alpha\not\sim_J u$ and $\alpha\sim_I u$ then $\sim_I$ is generated by $\sim_J$ and $\alpha\sim_I u$.}
}
\end{center}
is also true: apply point b) after exchanging $I$ and $J$ and after replacing $\varphi$
by $\varphi^{-1}=\varphi_{\alpha,w,-\tau}$.
\end{remark}
\vskip 5pt
If $\sim$ and $\sim'$ are homotopy relations,
we will say $\sim'$ is a \textbf{direct successor} of $\sim$ if there exist quivers with admissible relations $(Q,I)$ and
$(Q,J)$ presenting $A$ together with $\varphi_{\alpha,u,\tau}$ a transvection such that: $\sim_I=\sim$,
$\sim'=\sim_J$, $J=\varphi_{\alpha,u,\tau}(I)$, $\alpha\not\sim_I u$ and $\sim_J\ =\ <\,\sim_I\ ,\ u\sim_J \alpha>$ .
Notice that there may exist various $(Q,I)$, $(Q,J)$ and $\varphi_{\alpha,u,\tau}$ providing the same homotopy relations
and a direct successor relation between
them. These remarks will be taken into account in the following definition.
%
%
%
%
%
%
\section{Proof of Theorem \ref{THM}}
Assume $A$ is a finite dimensional basic and connected $k$-algebra with ordinary quiver $Q$.\\
\begin{e-definition}
If $Q$ has no oriented cycles, we define a quiver $\Gamma$ as follows:
\begin{enumerate}
\item[.] the vertices of $\Gamma$ are the homotopy relations of the admissible presentations of $A$ with quiver and relations,
\item[.] $\Gamma$ has an arrow $\sim \to \sim'$ if and only if $\sim'$ is
a direct successor of $\sim$.
\end{enumerate}
\end{e-definition}
\vskip 10pt
The author thanks Mariano Su\'arez-Alvarez for the following remark:
\vskip 5pt
\begin{remark}
A homotopy relation is determined by its restriction to the paths of
$Q$ with length at most the radical length of $A$, thus there are only
finitely many homotopy relations. This argument shows that $\Gamma$ is
finite. As a consequence, any vertex of $\Gamma$ is the target of a
(finite) oriented path in $\Gamma$ with source a source of $\Gamma$
(i.e. a vertex with no arrow arriving at it).
\end{remark}
\vskip 5pt
Moreover, using Lemma \ref{lem2}, the following result gives
additionnal properties of the quiver $\Gamma$.\\
\begin{e-proposition}
\label{prop1}
Assume $Q$ has no oriented cycles and let $m$ be the number of bypasses
in $Q$. Then
$\Gamma$ is connected and  has no oriented cycles. Any vertex of $\Gamma$ is the source
of at most $m$ arrows. The length of the oriented paths in $\Gamma$
is bounded by $m$.
\end{e-proposition}
\vskip 10pt
Suppose now that $\sim_{I_0}\to\ldots\to\sim_{I_n}$ is an oriented path in $\Gamma$.
For each $i$, the natural map $\pi_1(Q,I_i)\to \pi_1(Q,I_{i+1})$
induced by the identity map on walks is a well defined surjective group morphism. Thus there
is a surjective group morphism $\pi_1(Q,I_0)\twoheadrightarrow\pi_1(Q,I_n)$.
Consequently, it is natural to ask whether $\Gamma$ has a unique source.
 In case of a positive answer, this unique source gives rise to a canonical homotopy relation among all other homotopy relations of the presentations of the given algebra $A$.
The following Proposition answers this question. 
Notice that the preceding results did not use the hypotheses (written
before stating Theorem \ref{THM}) concerning $car(k)$ or
concerning the double bypasses.\\
\begin{e-proposition}
\label{prop2}
Assume $A$ satisfies the hypotheses written before stating Theorem \ref{THM},
then $\Gamma$
has a unique source.
\end{e-proposition}
\vskip 5pt
The proof of Theorem \ref{THM} is now straightforward: let $kQ/I_0\simeq A$ be a presentation
of $A$ such that $\sim_{I_0}$ is the unique source of $\Gamma$. For any other homotopy relation
$\sim_J$ (with $kQ/J\simeq A$), there exists a path $\sim_{I_0}\to \ldots \to \sim_J$ in $\Gamma$,
thus we have a surjective group morphism $\pi_1(Q,I_0)\twoheadrightarrow\pi_1(Q,J)$.
\vskip 5pt
\begin{remark}
If $m$ is the number of bypasses of $Q$, the unicity of the source of
 $\Gamma$ implies that $\Gamma$ has at most $1+m+m^2+\ldots+m^m$ vertices.
In particular, under the assumptions made before stating Theorem \ref{THM}, there are at most $1+m+m^2+\ldots+m^m$ isomorphism
classes of groups which can be the fundamental group of a presentation of $A$ with quiver and admissible relations.
\end{remark}
\begin{remark}
Let $Q=$
\psset{unit=10pt}
\begin{pspicture}(0,0)(8,3)
\rput(0,0){\rnode{A}{}}
\rput(2,2){\rnode{B}{}}
\rput(4,0){\rnode{C}{}}
\rput(6,2){\rnode{D}{}}
\rput(8,0){\rnode{E}{}}
\ncline{->}{A}{B}
\ncline{->}{A}{C}
\bput{0}{$a$}
\ncline{->}{B}{C}
\ncline{->}{C}{D}
\ncline{->}{C}{E}
\bput{0}{$\alpha$}
\ncline{->}{D}{E}
\end{pspicture}.
Notice that $Q$ has no double bypasses.
Let $u$ (resp. $v$) be the path \vskip 5pt parallel to $a$ (resp. $\alpha$), let $I_1=<\alpha a+vu,va+\alpha u>$
and set $A\simeq kQ/I_1$. Then $\pi_1(Q,I_1)=\mathbb{Z}/2$. Let $I_2=\varphi_{a,u,-1}\varphi_{\alpha,v,-1}(I_1)$,
hence $A\simeq kQ/I_2$. If $car(k)=0$, then $I_2=<\alpha a,va+\alpha u-2vu>$, $\pi_1(Q,I_2)=0$ and $\sim_{I_1}$ is the unique source of $\Gamma$.
Suppose now that $car(k)=2$. Then $I_2=<\alpha a,va+\alpha u>$ and $\sim_{I_1}$
and $\sim_{I_2}$ are both sources of $\Gamma$ whereas $\sim_{I_1}=\,<\sim_{I_2},\alpha a\sim_{I_1}vu>$ does not coincide
with $\sim_{I_2}$. Notice that we still have a surjection $\pi_1(Q,I_2)=\mathbb{Z}\twoheadrightarrow \mathbb{Z}/2=\pi_1(Q,I_1)$.
Notice also that one can build similar examples for any non zero value $p$ of $car(k)$ by taking for $Q$ a sequence of $p$ bypasses instead of $2$ bypasses only.
\end{remark}
\vskip 5pt
\begin{remark}
The details of our proof of Theorem \ref{THM} show that the assumption on $car(k)$ can be weakened. More precisely,
if $p=car(k)\neq 0$ and if the quiver $Q$ has less than $p$ bypasses, then Theorem \ref{THM} still holds.
\end{remark}
\vskip 5pt
The results presented here can be reformulated into results on Galois coverings (where the group
$\pi_1(A)$ corresponds to the universal cover of $A$ and which does not depend on the presentation of $A$).
This reformulation will be made in a subsequent paper.
%
%
%
%
%
%
\section*{Acknowledgements}
The author gratefully acknowledges Claude Cibils and Mar\'ia Julia Redondo for many helpful comments and for a careful reading
of the preliminary versions of this note.
%
%
%
%
%
%

\end{document}